\documentclass[11pt,onecolumn,journal]{article}
\usepackage[bottom=2.5cm,top=2.5cm,left=2.5cm,right=2.5cm]{geometry}

\usepackage{amsmath,amssymb,amsthm,graphicx,caption,cases}	
\usepackage[bottom]{footmisc}	
\usepackage{algorithmic}

\usepackage{abstract,amsfonts,amsbsy,amssymb,amstext,mathrsfs,latexsym,color,url,amsthm,xcolor,authblk,fontenc,bm}

\definecolor{NUSBlue}{RGB}{0,61,124} 
\definecolor{NUSOrange}{RGB}{239,124,0}
\usepackage[colorlinks=false,allbordercolors=NUSOrange]{hyperref}

\DeclareOldFontCommand{\bf}{\normalfont\bfseries}{\mathbf}

\newcommand{\bp}{\mathbf{p}}
\newcommand{\bz}{\mathbf{z}}

\newcommand{\bv}{\mathbf{v}}

\newcommand{\bb}{\mathbf{b}}

\newcommand{\bu}{\mathbf{u}}

\newcommand{\bs}{\mathbf{s}}

\newcommand{\bw}{\mathbf{w}}
\newcommand{\bx}{\mathbf{x}}
\newcommand{\by}{\mathbf{y}}

\newcommand{\R}{\mathbb R}

\definecolor{NUSBlue}{RGB}{0,61,124}   

\usepackage{tikz}
\usepackage{mathtools}

\def\nudge{.5}

\tikzset{axis/.style={ultra thin, Grey, -latex, shorten <=-\nudge cm, shorten >=-2*\nudge cm}}
\tikzset{line/.style={thick}}

\usepackage{textcomp}
\DeclareMathAlphabet\mathbfcal{OMS}{cmsy}{b}{n}
\usepackage{hyperref}
\usepackage[linesnumbered,vlined,ruled]{algorithm2e}
\usepackage{amsmath,amsfonts,amssymb,amsthm}
\theoremstyle{plain}
\newtheorem{thm}{Theorem}

\newtheoremstyle{cited}%
  {3pt}
  {3pt}
{\itshape}
  {}
  {\bfseries}
  {.}
  {.5em}
  {\thmname{#1} \thmnumber{#2} \thmnote{\normalfont#3}}
\theoremstyle{cited}
\newtheorem{citedthm}[thm]{Theorem}
\newtheorem{citedlem}[thm]{Lemma}

\newtheorem{citeddef}[thm]{Definition}
\newtheorem{citedprop}[thm]{Proposition}

\usepackage{graphicx,epstopdf,epsf,subcaption} 
\usepackage{cite}
\usepackage{tabularx}
\usepackage{array}
\usepackage{booktabs}
\usepackage{comment}
\usepackage{float}
\usepackage{caption}

\begin{document}
\renewcommand*{\Authsep}{, }
\renewcommand*{\Authand}{, }
\renewcommand*{\Authands}{, }
\renewcommand*{\Affilfont}{\normalsize}   
\setlength{\affilsep}{2em}   
\title{Convex optimization and block coordinate descents using direction  proximal point method}
\title{Directional proximal point method for convex optimization}
\date{}
\author[1]{Wen-Liang Hwang}
\author[2]{Chang-Wei Yueh}
\affil[1]{The Institute of Information Science, Academia  Sinica, Taipei, Taiwan.}
\affil[2]{The Department of Electrical Engineering, National Taiwan University, Taipei, Taiwan. }
\maketitle

\begin{abstract}
The use of proximal point operators for optimization can be computationally expensive when the dimensionality of a function (i.e., the number of variables) is high. In this study, we sought to reduce the cost of calculating proximal point operators by developing a directional operator in which the proximal regularization of a function along a specific direction is penalized. We used this operator in a novel approach to optimization, referred to as the directional proximal point method (Direction PPM). When using  Direction PPM, the key to achieving convergence is the selection of direction sequences for directional proximal point operators. In this paper, we present the conditions/assumptions by which to derive directions capable of achieving global convergence for convex functions. Considered a light version of PPM, Direction PPM uses scalar optimization to derive a stable step-size via a direction envelope function and an auxiliary method to derive a direction sequence that satisfies the assumptions. This  makes Direction PPM adaptable to a larger class of functions. Through applications to differentiable convex functions, we demonstrate that negative gradient directions at the current iterates could conceivably be used to achieve this end. We provide experimental results to illustrate the efficacy of Direction PPM in practice. 
\end{abstract}

\section{Introduction}

The proximal point method (PPM) is a state-of-the-art method used in optimization \cite{bolte2014proximal}, signal processing \cite{becker2011nesta}, distributed convex optimization \cite{BoydADMMsurvey}, and machine learning \cite{zhang2018ista}. It involves the iterative application of convex optimization updates to a fixed point of the proximal point operator, which represents a minimum for the convex function. 
The proximal point operator is a theoretically sound and elegant approach to the construction of optimization algorithms. 
Some applications involve simple or structured functions that allow for the efficient assessment of their proximal operators either via closed-form expressions or simple linear time algorithms \cite[Chapter 6]{parikh2014proximal}. However, for most functions, the calculation of their proximal point operators imposes an enormous computational burden \cite[Section 3.2]{parikh2014block}.
The effects are particularly obvious when the dimensionality (i.e., the number of variables) is high. 
For example, the proximal point operator for proper, closed, and convex function $f: \R^n \rightarrow \R \cup \{\infty\}$ is defined as follows:
\begin{align} \label{proxop}
\text{prox}_{tf}(\bx) = \arg \min_{\bu \in \R^n} f(\bu) + \frac{1}{2t} \| \bu - \bx \|^2.
\end{align}
The proximity operator can be characterized using the sub-differential of $f$ via
\begin{align}\label{PSDcharacterization}
x-\text{prox}_{tf}(x) =  t\partial f(\text{prox}_{tf}(x)),
\end{align}
where the sub-differential of $f$ is the set-valued operator on $\R^n$ given by
\begin{align*}
\partial f(x)= \{ u\in\R^n \colon \,\langle y-x,u\rangle + f(x) \leq f(y)
\text{ for all }  y\in \R^n \}.
\end{align*}

This operation searches for $t\partial f(\text{prox}_{tf}(x))$ in $\R^n$, which is the vector from the following iterate $\text{prox}_{tf}(\bx)$ to the current iterate $\bx$. This search can be costly when $n$ is large, because optimizing the objective function for proximal regularization necessitates searching for both the direction and step-size in search space $\R^n$.

We thus suggest reducing the cost of calculating the proximal point operator by splitting the search into two separate tasks. 
By determining the direction from one iterate to the next, we can transform the search for a suitable step-size into a scalar optimization problem. This approach to separation is well-suited to large-scale systems in which the increase in optimization cost exceeds what is tractable using proximal point operators \cite{asl2020analysis,fletcher2005barzilai}.
We define the direction envelope of convex function $f$ with respect to a direction $\bar \bp$ (unit-norm direction) and parameter $t > 0$ as follows:
\begin{align} \label{dppmwform}
f_t^{\bar \bp}(\bx) = \inf_{w \geq 0}\frac{1}{2t} \|w\|^2 + f(w \bar \bp + \bx).
\end{align}
The infimum in (\ref{dppmwform}) is achieved at a unique point because the quadratic regularization term is supercoercive on $w \geq 0$. Accordingly, we define the directional proximal point operator $\text{prox}_{tf}^{\bar \bp}: \R^n \rightarrow \R^n$ with the iterate updated to
\begin{align} \label{dppmdefinition}
\text{prox}_{tf}^{\bar \bp}(\bx) = \bx + \bar \bp [\arg f_t^{\bar \bp} (\bx)].
\end{align}
A comparison of (\ref{proxop}) with (\ref{dppmdefinition}) shows the profound difference between PPM and the proposed approach (referred to as Direction PPM). In Direction PPM, the search direction is derived using the current iterate, and the optimal step-size is derived using the following iterate by solving (\ref{dppmwform}), while both the search direction and step-size are derived from the following iterate via PPM. The sequence of directions $\{\bar \bp_k\}$ and parameters $\{t_k\}$ in Direction PPM directs the algorithm toward convergence. We demonstrate that the value of $t$ is relevant to the optimal step-size but irrelevant to convergence. If convergence was the sole concern, then the value of $t$ could be set 
as a constant. Note however that the direction of the sequence is crucial to the convergence of the algorithm. In this paper, we present two conditions sufficient for direction sequences to reach convergence: sub-gradient-relatedness and target-relatedness. The assumption of sub-gradient-relatedness prevents the direction sequence from tending toward orthogonality to any sub-gradient at a limit point, and the assumption of target-relatedness ensures that the search direction has a sub-component targeting the critical points of a function. We demonstrate that a descent direction at an iterate can be used as a Direction PPM search direction by which to achieve sub-sequence convergence. Furthermore, mild conditions on the function should be sufficient to achieve convergence of the entire sequence from the sequence of negative sub-gradient descent directions. For differentiable convex functions, negative gradient directions at current iterates can be used in Direction PPM to achieve convergence, wherein the convergence rate (i.e., the order of updates required to obtain sub-optimal solutions) matches that of PPM. Nonetheless,  
executing an update for PPM is generally far more computationally expensive than for Direction PPM. 
Direction PPM can be used for the optimization of convex functions, or used in conjunction with  alternating optimization for non-convex functions. This makes it possible to leverage the computational efficiency and stability of the step-size selection for first-order optimization methods.

In Section \ref{backgrounds}, we compare Direction PPM with other first-order optimization methods that are applicable to convex functions.  
In Section \ref{sec:DPPM}, we present elementary properties of direction envelope functions and derive the optimal step-size for an update. In Section~\ref{sec:generalDPPM}, we present the conditions for the direction sequences of Direction PPM, which are related to the algorithm convergence. In Section~\ref{sec:experiments}, we present experimental results showing that Direction PPM works in practice as theory asserts, and is competitive against existing algorithms. In Section \ref{cons}, we present concluding remarks.

 {\bf{Notation:}} \\
\noindent$\bullet$ $\|\bx\|$ denotes the 2-norm of $\bx$ and $\bar \bp$ denotes the unit 2-norm vector of $\bp$. \\
\noindent$\bullet$ Boldfaced letters are vectors. \\
\noindent$\bullet$ $\Gamma$ denotes the class of proper, closed, and convex function; and $ \{\bx| f(\bx) \leq f(\bx^0)\}$ is a bounded set for a given $\bx^0$.

\section{Related works} \label{backgrounds}

Standard first-order methods for unconstrained optimization problems rely on ``iterative descent", which works as follows: Start at point $\bx^0$, and successively generate iterates $\bx^1$, $\bx^2$, $\cdots$, such that function $f$  decreases in each iteration. The aim is to decrease $f$ to its minimum.  In many cases, this can be achieved by selecting a step-size in the search direction based on the properties of $f$ \cite{nocedal2006numerical, Ber08}. Direction PPM is a first-order optimization method. In the following, we compare Direction PPM with other first-order optimization methods applicable to convex functions.

\subsection{Differentiable convex functions}

A popular approach to differentiable $f  \in \Gamma$ is the gradient descent method, where the selection of step-size is based on the backtracking procedure, in which the step-size decreases successively  until the Armijo condition is met. This approach does not generate a sequence capable of converging to a single critical point. In contrast, Direction PPM uses the negative gradient direction to ensure that the entire sequence of iterates converges to a single critical point of $f$ (as indicated in Theorem \ref{dppmfejer}).

The Barzilai and Borwei method is appropriate for large-scale optimization problems  \cite{barzilai1988two} as long as efficiency is the sole consideration. Note that this method requires only ${\cal O}(n)$ floating point operations and a gradient evaluation to update in $\R^n$. The search direction of the method is always along the direction of the negative gradient. The step-size is not a conventional choice for a back-tracking line search, as it cannot guarantee descent in every update. The main challenge of the Barzilai and Borwei method is stabilizing step-size selection to avoid long step-sizes, which could lead to divergence from the optimal solution  \cite{burdakov2019stabilized}. 
Direction PPM derives step-size by solving a scalar optimization problem. For differentiable functions, Direction PPM can converge to a single critical point as long as the negative gradients of the respective iterates are the designated search directions. The PPM is another option for large-scale problems, due to its stability in deriving solutions; however, for most functions, optimizations over $\R^n$ incur excessively high computational costs, as the calculation of each iterate involves the gradient ($\nabla f(\text{prox}_{tf}(x))$) of the following iterate.  Direction PPM achieves convergence for the entire sequence with a similar convergence rate to that of the PPM but with far lower computational costs, due to its use of scalar optimizations and gradients ($\nabla f(x)$) of the current iterates.

\subsection{Non-differentiable convex functions}

In the following, we compare sub-gradient and proximal point methods for non-differentiable convex functions. 

The sub-gradient method does not necessarily descend after each update, due to the fact that it always updates along a negative sub-gradient direction. This method cannot converge to critical points of $f$ using a constant step-size. However, for step-size sequence $\{t_i = \frac{1}{i+1}\}$ (satisfying $\sum_i t_i = \infty$ and $\sum_i t_i^2 < \infty$), this method achieves global convergence to a single critical point of a convex function. Note that this step-size sequence is not easy to use in  practice because it rapidly drops to zero. This has prompted the development of variants of the step-size selection strategy to deal with optimizations that involve multiple iterations with variable updates \cite{bottou2018optimization}. Note also that stopping the algorithm can be problematic, as it does not always descend.


Unlike the sub-gradient method, Direction PPM is a descent algorithm (Lemma \ref{descentstep}).
The use of sub-gradient directions that are also descent directions enables global convergence of the entire sequence of iterates to a single critical point of a convex function (Theorem \ref{dppmfejer}). This makes it possible to use $\|f(\bx^{k+1}) - f(\bx^k)\| \leq \epsilon$ as a stopping condition.
If the sole consideration was the minimization of the computational complexity involved in reaching an $\epsilon$-suboptimal solution, then Direction PPM would be superior to the sub-gradient method. The number of updates when using the sub-gradient method is in the order of ${\cal O)}(\frac{1}{\epsilon^2})$, whereas the number of updates when using Direction PPM is ${\cal O}(\frac{1}{\epsilon})$, even without acceleration. 

We can consider Direction PPM a light version of PPM. Direction PPM uses scalar optimization to derive a stable step-size via a direction envelope function as well as an auxiliary method to derive a direction sequence that satisfies sub-gradient-related and/or target-related assumptions.  
Generally, for functions in which it is easy to obtain one descent direction in each update, Direction PPM has an edge over PPM. 
If a descent direction at a given point is used as the search direction, then Direction PPM can converge to critical points (Proposition \ref{dppmlimitpoints}). Furthermore, if a sub-gradient at a given point is also a descent direction at that point, then the sub-gradient can be used as the search direction to achieve convergence to a single critical point of a function (Theorem \ref{dppmfejer}).

\section{Characterization of updates in Direction PPM}  \label{sec:DPPM}

In the following, we demonstrate that direction envelope function $f_t^{\bar \bp}(\bx)$ at search direction $\bar \bp$ has the same optimal value over $\bx \in 
R^n$, regardless of the value of $t$.

\begin{citedlem} \label{samemin}
For a given $\bar \bp$ and $f \in \Gamma$.  \\
(i) $f(\bx) \geq f_t^{\bar \bp}(\bx) \geq f_{t^+}^{\bar \bp}(\bx) \geq \min_{\bx} f(\bx)$ with $t^+ \in [t, \infty)$. \\
(ii) $\min_{\bx} f_t^{\bar \bp} (\bx) = \min_{\bx} f(\bx)$ for any $t > 0$.
\end{citedlem}
\proof
(i) Clearly, $f_{t}^{\bar \bp}(\bx) \leq f(\bx)$ (as is the case of $f(\bx + w \bar \bp)$ with $w = 0$). $f_t^{\bar \bp}(\bx) = f(\bx + w_t \bar \bp) + \frac{1}{2t}\| w_t\|^2 \geq f(\bx + w_t \bar \bp) + \frac{1}{2t^+}\| w_t\|^2 \geq  f_{t^+}^{\bar \bp}(\bx)$. 

(ii) From (i), 
we have $\min_{\bx} f^{\bar \bp}_t(\bx) \leq \min_{\bx} f(\bx)$ and $\min_{\bx} f_t^{\bar \bp}(\bx) \geq \min_{\bx} f(\bx)$.
\qed

The directional proximal point operator for a given search direction is updated in Eq. (\ref{dppmdefinition}), and the step-size along that direction is derived by solving scalar optimization problem (\ref{dppmwform}).

\subsection{Scalar optimization}

For direction $\bar \bp$ ($\| \bar \bp\| = 1$) at $\bx$ where parameter $t > 0$, we consider the Moreau envelope function\footnote{$\min_{\bu} \frac{1}{2t} \| \bu - \bx \|^2 + f(\bu)$} with constraints $\bu = \bx + w \bar \bp$ and step-size $w \geq 0$:
\begin{align} \label{dppmoptfor}
\begin{cases}
\min_{w,\bu} \frac{1}{2t} \| \bu - \bx\|^2 + f(\bu) \\
\bu - \bx = w\bar \bp \\
w \geq 0.
\end{cases}
\end{align}
Substituting the constraint into the objective yields the direction envelope function $f_t^{\bar \bp}(\bx)$ (\ref{dppmwform}) for Direction PPM.
Suppose that $f \in \Gamma$; then, $\frac{1}{2t} \|w\|^2 + f(w \bar \bp + \bx)$ is a strictly convex function of $w$. Thus, the solution to Eq. (\ref{dppmwform}) is unique and can be expressed as follows:
\begin{align} \label{dppmwform1}
w^* = \arg \min_{w \geq 0}\frac{1}{2t} \|w\|^2 + f(w \bar \bp + \bx).
\end{align}
The explicit form of $w^*$ is presented in the following lemma. Here, we need the sub-gradient upper bound $M$ of $f $, where $M \geq \| s(\bx) \|$ for all $s(\bx) \in \partial f(\bx)$ and all $\bx$.
\begin{citedlem}\label{dppmunique}
Let $M < \infty$ denote the sub-gradient upper bound for $f \in \Gamma$. 
Then, the solution of (\ref{dppmwform1}) is  
\begin{align}  \label{stepsizeval}
w^* = 
\begin{cases}
- t \bar \bp^\top \bs(\bu) > 0 \text{ if $\bar \bp^\top \bs(\bu)  < 0$} \\
0  \text{ otherwise}.
\end{cases}
\end{align}
 \end{citedlem}
\proof
In accordance with variational inequality \cite{kinderlehrer2000introduction}, $w^*$ is a solution to (\ref{dppmwform}) if and only if there exists a sub-gradient $\bs(\bu)$ of $f$ at $\bu$, such that
\begin{align} \label{dppmvarI}
(\frac{1}{t} w^* + \bar \bp^\top \bs(\bu)) (w - w^*) \geq 0 \text{ for all $w \geq 0$}.
\end{align}
This equation is equivalent to    
\begin{align} \label{dppmvi}
\begin{cases}
\frac{1}{t} w^* + \bar \bp^\top \bs(\bu) = 0 \text{ if $w^* > 0$} \\
 \bar \bp^\top \bs(\bu) \geq 0 \text{ if $w^* = 0$}.
\end{cases}
\end{align}
Therefore, we can assert (\ref{stepsizeval}).
The range of $w^*$ is derived based on 
\begin{align}\label{tparameter}
0 \leq w^*  & \leq t | \bar \bp^\top \bs(\bu) |  \leq t M.
\end{align}
\qed

Equation (\ref{tparameter1}) below shows that the global bound involving $M$ can be improved to a locally attainable bound. From the convexity of a function, it can be easily derived that $w^* = 0$ for any parameters $(t, \bar \bp)$. 
As shown in the following, updating the directional proximal point operators using Eq. (\ref{dppmdefinition}) does not increase function value in any direction, regardless of whether a convex function is differentiable or non-differentiable.

\begin{citedlem} \label{descentstep}
Suppose that $f \in \Gamma$ and that $M$ is the sub-gradient upper bound of $f$. Let
\begin{align} \label{dppmupdatek0}
\bu = \bx - t [\bs(\bu)^\top \bar \bp] \bar \bp = \bx + w^* \bar \bp,
\end{align}
where $w^*\in [0, tM]$ is the solution to (\ref{dppmwform1}) and $s(\bu) \in \partial f(\bu)$. 
Then, we have 
\begin{align}   \label{dppmdescent0}
f(\bu) & \leq f(\bx) - t | \bs(\bu)^\top \bar \bp|^2.
\end{align}
\end{citedlem}
\proof
Because $f$ is convex, we obtain 
\begin{align*}
f(\bz)  \geq f(\bu) + \bs(\bu)^\top (\bz - \bu) 
 = f(\bu)  + \bs(\bu)^\top  (\bz - \bx + t(\bs(\bu)^\top \bar \bp) \bar \bp).
\end{align*}
Re-writing the above yields 
\begin{align*}
f(\bu) & \leq  f(\bz)  + \bs(\bu)^\top (\bx  - \bz) - t | \bs(\bu)^\top \bar \bp|^2.
\end{align*}

Allowing $\bz = \bx$ gives us 
\begin{align} \label{dppmsufficient}
f(\bu) \leq f(\bx) - t | \bs(\bu)^\top \bar \bp|^2.
\end{align}
\qed

Numerous algorithms have been developed for scale optimization \cite{press1988numerical}. In the following, we review those that are easily implemented for optimal values within a closed interval. If $f$ is differentiable, it is possible to use the bisection method to derive the unique root of 
\begin{align*} 
0= w + t \bar \bp^\top \nabla f(\bx + w \bar \bp).
\end{align*}
If $f$ is not differentiable, the golden section search method \cite{avriel1968golden} can be used to find the unique minimum of Eq. (\ref{dppmwform}). The fact that the golden search method does not use any gradient evaluations makes it ideally suited to situations in which the gradient of a function cannot be efficiently or accurately derived as well as to situations n which the function is not differentiable. 
In the golden section method, the minimum is bracketed when there is a triplet  $(a, b, c)$ with $a < b < c$, such that $f(\bx + b \bar \bp)$ is larger than $f(\bx + a \bar \bp)$ and $f(\bx + c \bar \bp)$. This method involves selecting a new point $z$, either between $a$ and $b$ or between $b$ and $c$. If we make the latter choice, then we evaluate $f(\bx + z \bar \bp)$. If $f(\bx + b \bar \bp) > f(\bx + z \bar \bp)$, then the new bracketed triplet is $(a, b, z)$. If $f(\bx + b \bar \bp) < f(\bx + z \bar \bp)$, then the new bracketed triplet is $(b, z, c)$.  This is repeated until the bracket corresponding to the minimum is sufficiently small.

\subsection{Descent directions of Direction PPM}

We now present some useful properties of Direction PPM.

\begin{citedlem}\label{dppmmonoprop}
Suppose that $f \in \Gamma$. For $w_1, w_2 \in \R_{+}$, we let $\bu_1 = \bx + w_1 \bar \bp$ and $\bu_2 = \bx + w_2 \bar \bp$. Denote $\bs(\bz)$ as the sub-gradient of $f$, at $\bz$. This gives us the following: \\
(i) \begin{align} \label{dppmmonoprop1}
[\bs(\bu_1)^\top \bar \bp - \bs(\bu_2)^\top \bar \bp]\bar \bp^\top (\bu_1 - \bu_2) \geq 0.
\end{align}
(ii) $\bar \bp^\top \bs(\bx + w \bar \bp)$ is a monotonically non-decreasing function of $w \in \R_{+}$. For $w \in [0, w^*]$ (the optimal step-size defined in (\ref{dppmwform1})), $\bar \bp^\top \bs(\bx + w \bar \bp) \leq 0$ and therefore
\begin{align} \label{monodir}
| \bar \bp^\top \bs(\bx) | \geq |\bar \bp^\top \bs(\bx + w \bar \bp)|.
\end{align}
(iii) $f(\bx + w \bar \bp) + \frac{1}{2t}\|w\|^2$ is a decreasing function on $[\bx, \bx + w^* \bar \bp]$ with $w^* > 0$ and $w \in [0, w^*]$.
\end{citedlem}
\proof
From $\bu = \bx + w \bar \bp$, 
we obtain $w = \bar \bp^\top (\bu - \bx)$. 
\begin{align}  \label{dppmmonoprop0}
[\bs(\bx) - \bs(\bu)]^\top (\bx - \bu)  & =  - w [\bs(\bx) - \bs(\bu)]^\top \bar \bp  \nonumber \\
& = [(\bs(\bx)^\top \bar \bp) - (\bs(\bu)^\top \bar \bp)]\bar \bp^\top (\bx - \bu) \geq 0.
\end{align} 
The last inequality is due to the monotonicity property of the sub-gradients of convex functions.
Without a loss of generality, we can suppose that $w_2  > w_1$. Substituting $\bu_1 = \bx + w_1 \bar \bp$ and $\bu_2 =  \bu_1 + (w_2 - w_1) \bar \bp$ respectively for $\bx$ and $\bu$ into (\ref{dppmmonoprop0}) yields (\ref{dppmmonoprop1}).

(ii) If substituting $\bu_1 - \bu_2 = (w_1 - w_2) \bar \bp$ into (\ref{dppmmonoprop1}), we obtain  
\begin{align*}
(w_1 - w_2)  (\bs(\bu_1) - \bs(\bu_2))^\top \bar \bp \geq 0.
\end{align*}
Because $w_2  > w_1$, we obtain $\bs(\bu_1)^\top \bar \bp \leq  \bs(\bu_2)^\top \bar \bp$. Therefore, $\bar \bp^\top \bs(\bx + w \bar \bp)$ is a monotonically non-decreasing function of $w \in \R_{+}$.  Because $\bar \bp^\top \bs(\bx + w^* \bar \bp) \leq 0$, we obtain (\ref{monodir}) for $w \in [0, w^*]$.

(iii) Taking the derivative of $f(\bx + w \bar \bp) + \frac{1}{2t}\|w\|^2$ with respect to $w$, we obtain $\frac{w}{t}+ \bar \bp^\top \bs(\bx + w \bar \bp)$. From part (ii), $\bar \bp^\top \bs(\bx + w \bar \bp)$ is a  monotonically non-decreasing function of $w$. Therefore, $\frac{w}{t}+ \bar \bp^\top \bs(\bx + w \bar \bp)$ increases from $\bx$ to $\bx + w \bar \bp$. Because envelope function $f^{\bar \bp}_t(\bx)$ obtains the minimum value of (\ref{dppmwform}) at $\bx + w^* \bar \bp$, $f(\bx + w \bar \bp) + \frac{1}{2t}\|w\|^2$ decreases from $\bx$ to $\bx + w^* \bar \bp$.

\qed

In accordance with  (\ref{monodir}), the upper bound in (\ref{tparameter}) can be refined to obtain
\begin{align}\label{tparameter1}
0 \leq w^*  & \leq t | \bar \bp^\top \bs(\bu) |  \leq t | \bar \bp^\top \bs(\bx) | \leq t |f^{'}(\bx, \bar \bp)|,
\end{align}
where $f^{'}(\bx, \bar \bp)$ is the directional derivative of $f$ at $\bx$ in the direction $\bar \bp$.
Lemma \ref{descentstep} indicates that Direction PP would not increase function values, regardless of the search direction. Hence, we must characterize the descent directions that lead to a strict decrease in function values (i.e., $f(\bu) < f(\bx)$ with $\bu =  \bx + w^* \bar \bp$).  Although the notion of descent direction is well-defined for a differentiable function (any direction $\bp$ that leads to $\nabla f(\bx)^\top \bar \bp < 0$), the notion is not defined for a non-differentiable function. For non-differentiable function $f$, the following definition is a sufficient condition for $f(\bu) < f(\bx)$ based on $
f(\bx) \geq f(\bu) + \bs(\bu)^\top (\bx - \bu):
$
\begin{citeddef} \label{directdef}
Suppose that there is a sub-gradient $\bs(\bu) \in \partial f(\bu)$ and $\varepsilon >0$, such that $\bs(\bu)^\top (\bx - \bu) = -\bs(\bu)^\top \bp \geq \varepsilon$, then $-\bs(\bu)^\top \bar \bp \geq \frac{\varepsilon}{\|\bp \|}$. Hence, $\bar \bp$ is referred to as the DPPM-descent direction of parameter $\varepsilon$ for function $f$ at $\bx$. For brevity, we refer to $\bar \bp$ as a DPPM-descent direction or simply a descent direction (when this will not cause confusion) to indicate that $\bar \bp$ is a DPPM-descent direction for some $\varepsilon > 0$. 
\end{citeddef}

Applying the above definition to differentiable function $f \in \Gamma$ with $\bp = -\nabla f(\bx)$ with $\| \nabla f(\bx) \|^2 > \varepsilon$ at non-minimizer $\bx$ and using the fact that $\nabla f$ is a continuous function, we can derive a sufficiently small $v > 0$, such that $\bu = \bx - v \nabla f(\bx)$ and $\nabla f(\bu)^\top \nabla f(\bx) \geq \varepsilon$. Hence, $-\nabla f(\bx)$ is the DPPM-descent direction at $\bx$.


Clearly, a DPPM-descent direction of a function at a point is also a descent direction at the point of the function.
Suppose that $\bar \bp$ is a DPPM-descent direction. Using (\ref{stepsizeval}) and $\bu = \bx + w^* \bar \bp$, we can derive a non-zero optimal step-size $w^*$ of (\ref{dppmwform1}) with the following range of $t$: 
\begin{align} \label{tlowerbound}
\frac{L_{\bx}}{|f^{'}(\bx, \bar \bp)|} \geq  t  \geq   \frac{\bar \bp^\top (\bu - \bx)}{-\bs(\bu)^\top \bar \bp} = \frac{\| \bu -\bx\|}{-\bs(\bu)^\top \bar \bp},
\end{align}
where $L_{\bx}$ is the diameter (the largest segment) of the sub-level set $\{ \bz | f(\bz) \leq f(\bx) \}$ and the first inequality is derived from (\ref{tparameter1}) and $L_{\bx} \geq w^*$.
From the given fact that $\bar \bp$ is a DPPM-descent direction of $f \in \Gamma$ at $\bx$, we can deduce that $w^* > 0$ when $t$ is sufficiently large. The following equivalent statement can be made that $\bar \bp$ is not a descent direction of the Direction PPM: if $w^* = 0$ even for $t$ in (\ref{tlowerbound}), then $\bar \bp$ is not a descent direction of $f$ at $\bx$.

\section{Convergence analysis}\label{sec:generalDPPM}

Characterization of the descent directions at a given point for Direction PPM is insufficient to achieve convergence. Various assumptions pertaining to the search direction sequences are also necessary.

\subsection{Assumption of sub-gradient-relatedness}

To achieve convergence, we assume that the DPPM-descent directions for convex functions are 
sub-gradient-related. 
\begin{citeddef} \label{defgradientrelated}
Denote the descent direction of $f \in \Gamma$ at $\bx^k$ as $\bar \bp_k$. Let $\bs(\bz)$ be any sub-gradient of $f$ at $\bz$. The sequence $\{\bar \bp_k\}$ is sub-gradient-related, if for any subsequence $\{\bx^k\}_{k \in K}$ converging to a ``non-critical" point of $f$, the corresponding $\{\bar \bp_k\}_{k \in K}$ satisfies 
\begin{align} \label{gradient-related}
\lim_{k \rightarrow \infty} \sup_{k \in K}\bar \bp_k^\top \bs(\bx^{k+1}) < 0.
\end{align}
\end{citeddef}
This is equivalent to the following: $\{\bar \bp_k\}_{k \in K}$ converges to a critical point of $f$ when the descent sequence satisfies $\lim_{k \rightarrow \infty} \sup_{k \in K}\bar \bp_k^\top \bs(\bx^{k+1}) \rightarrow 0$.

\begin{citedprop} \label{dppmlimitpoints}
Suppose that $f \in \Gamma$.  Denote as $\{ \bx^k \}$ the sequences of iterates generated by Direction PPM using $\{(\bar \bp_k, t_k)\}_k$ at which $\bar \bp_k$ is a DPPM-descent direction with sub-gradient-relatedness, $t_k$ satisfies (\ref{tlowerbound}), and $\lim\inf_k t_k > 0$. This gives us the following: \\
(i) Asymptotic regularity (i.e. $\lim_{k \rightarrow \infty} \| \bx^{k+1}- \bx^k \| \rightarrow 0$) can be obtained with the entire sequence. \\
(ii) Every limit point of $\{ \bx^k \}$ is a critical point of $f$. \\
(iii) [Gradient-relatedness] For differentiable $f$, any subsequence $\{\bx^k\}_{k \in K}$ with corresponding descent direction sequence $\{\bp_k\}_{k \in K}$ that results in convergence to a non-critical point of $f$ must satisfy the following:
\begin{align} \label{subrelated1}
\lim_{k\rightarrow \infty} \sup_{k \in K} \nabla f(\bx^k)^\top \bp_{k} < 0.
\end{align}
\end{citedprop}

\proof

In accordance with Lemma \ref{descentstep} and (\ref{tlowerbound}), which gives a bound of $t_i$ for descent direction $\bar \bp_i$, we obtain
\begin{align}   \label{dppmdescent}
f(\bx^{i+1})  \leq f(\bx^{i}) - t_i| \bs(\bx^{i+1})^\top \bar \bp_i|^2 
 = f(\bx^{i}) -  w_i | \bs(\bx^{i+1})^\top \bar \bp_i|  
\end{align}
Summing the above results with $i=0, \cdots, k-1$, we obtain
\begin{align}
f(\bx^{k}) - f(\bx^0) &  \leq  -  \sum_{i=0}^{k-1} w_i | \bs(\bx^{i+1})^\top \bar \bp_i|   < \infty.
\end{align}
The right-hand side of the above equation is derived from the fact that $f$ is bounded from below. This implies that $w_k | \bs(\bx^{k+1})^\top \bar \bp_k|  \rightarrow 0$ as $k \rightarrow \infty$. Because $w_k = - t_k \bs(\bx^{k+1})^\top \bar \bp_k$, we obtain (i) with both
\begin{align} \label{sufficientcon}
w_k = \| \bx^{k+1} - \bx^k \| \rightarrow 0 \text{ and }  \bs(\bx^{k+1})^\top \bar \bp_k \rightarrow 0. 
\end{align}
(ii) Due to the assumption of sub-gradient-relatedness and (\ref{sufficientcon}), we conclude that any convergence subsequence $\{\bx^k\}_{k \in K}$ shall converge to critical points of $f$.  Moreover, (iii) and (\ref{subrelated1}) are a consequence of (ii) based on the fact that $\nabla f$ is a continuous function and (i) by replacing $\bx^{k+1}$ with $\bx^k$ in (\ref{gradient-related}).

\qed

For differentiable $f$, the ``non-orthogonality" condition (\ref{subrelated1}) is referred to as the gradient-related assumption \cite[(1.13)]{Ber08}, by which the sequence of iterates leading to convergence to the critical points of $f$ can be derived.

\subsection{Assumption of target-relatedness}

Imposing more than the assumption of sub-gradient-relatedness on search direction sequences allows us to establish the Fejer monotone of iterates, from which it is possible to demonstrate that the entire sequence of Direction PPM iterates can converge to a single critical point of $f$.


Let $C(\bx^0)$ denote the non-empty set of critical points to which iterates $\{\bx^k\}$ at initial $\bx^0$ converge and let $\bar \bp^*$ (i.e., $\overline{\bx^* - \bx})$ denote the direction from $\bx$ to $\bx^* \in C(\bx^0)$.
Suppose that $\bar \bp$ is the descent direction at $\bx \in \{\bx^k\}$. Then, $\bu = \bx - t [\bs(\bu)^\top \bar \bp] \bar \bp$ is the update.
If $\bs(\bu) \in \partial f(\bu)$ is in line with $\bar \bp$, then Direction PPM is PPM. Thus, in the following, we examine the situation where $\bs(\bu)$ and $\bar \bp$ are not parallel. Let $Q$ denote the two-dimensional plane spanned by $\bar \bp$ and $\bs(\bu)$. Denote $\bv = \mathbb P_Q(\bar \bp^*)$, which is the orthogonal projection of $\bar \bp^*$ to $Q$. Sub-gradient $\bs(\bu)$ can be expressed as follows:
\begin{align} \label{gradientdecomp}
\bs(\bu)= \alpha \bar \bp + \beta  \bv =  \alpha \bar \bp + \beta \bar \bv \| \bv\|.
\end{align}
Applying the inner product with $\bar \bp$ to both sides of the above equation yields
\begin{align}  \label{betaval}
0 > \bar \bp^\top \bs(\bu)= \alpha + \beta (\bar \bp^\top \bar \bv) \| \bv\|,
\end{align}
where the inequality can be attributed to the fact that $\bar \bp$ is a descent direction at $\bx$. 
We can also assume that
\begin{align} \label{directionassign}
 \bar \bp^\top \bar \bv> 0 \text{ and } \bs(\bu)^\top \bar\bp \geq \alpha 
 \end{align}
to obtain
\begin{align} \label{cosine}
\alpha < 0 \text{ and } 0 \leq \beta \leq  \frac{-\alpha}{(\bar \bp^\top \bar \bv) \| \bv\|} \text{ and }
\bar \bp =\frac{-1}{|\alpha|} \bs(\bu) + \frac{\| \bv\|\beta}{|\alpha|} \bar \bv.
\end{align}
Search direction $\bar \bp$ is in the quadrant enclosed by rays $-\bs(\bu)$ and $\bar \bv$ associated respectively with positive coefficients $\frac{1}{|\alpha|}$ and $\frac{\| \bv\|\beta}{|\alpha|}$.
We can deduce that 
\begin{align} \label{equalproj}
\bar \bp^\top \bar \bp^*= \bar \bp^\top (\mathbb P_Q(\bar \bp^*) + \mathbb P_{Q^\perp} (\bar \bp^*)) = \bar \bp^\top \mathbb P_Q(\bar \bp^*) =  \bar \bp^\top \bar \bv > 0,
\end{align} 
where the last inequality is due to (\ref{directionassign}), and
$\bar \bp$ has a sub-component pointing toward critical point $\bx^*$. Based on (\ref{stepsizeval}), $\bs(\bu)^\top \bar\bp \geq \alpha$ implies that $\bs(\bu)^\top \bar \bp = \frac{-w^*}{t}  \geq \alpha$;  therefore, $w^* \leq (-\alpha) t$. 
In accordance with Lemma \ref{dppmunique}, we can set $|\alpha| = M$.

The entire sequence of iterates converges to a single critical point of $f$ provided that $\beta \geq 0$ holds. Note that at $\beta = 0$, Direction PPM is PPM. Thus, we make an additional assumption pertaining to the descent direction in order to obtain $\beta \geq 0$:

\begin{citeddef} \label{targetdef}
DPPM-descent direction $\bar \bp$ at $\bx$ is target-related if the direction satisfies $\bar \bp^\top \bar \bp^* > 0$ (i.e., $\bar \bp^\top \bar \bv > 0$), where $\bp^* = \bx^* - \bx$ and $\bx^* \in C(\bx^0)$. 
\end{citeddef}

Because $f(\bx^*)  \geq f(\bx) + \bs(\bx)^\top (\bx^* - \bx) = f(\bx) + \bs(\bx)^\top \bp^*$ with $\bs(\bx) \in \partial f(\bx)$, any search direction $\bp = -\bs(\bx)$ with $\bs(\bx) \in \partial f(\bx)$ at a non-minimizer satisfies $\bp^\top\bar \bp^* > 0$. In conjunction with Definition \ref{targetdef}, we can derive that a DPPM-descent direction that is also a negative sub-gradient direction at a point is the target-related direction at the point. We note that there have been some numerical attempts to derive the sub-differentiable set at any non-minimizer for convex functions \cite{xia1992finding}. 

A Fejer monotone sequence implies that each iterate in the sequence is not strictly farther than its predecessor from any critical points of $f$, which means that the norm sequence $\{\| \bx^{k} - \bx^* \|\}$ converges for all $\bx^*$ in the set of critical points \cite{bauschke2011convex}.  In the following, we show that the entire sequence of $\{\bx^k\}$ converges to a single critical point of $f$. We achieve this by showing that the sequence of iterates of Direction PPM satisfies the Fejer monotone with respect to $C(\bx^0)$; i.e.,  for all $k$ and all $\bx^* \in C(\bx^0)$, $\| \bx^{k+1} - \bx^* \| \leq \|\bx^k - \bx^*\|$.

\begin{citedthm} \label{dppmfejer} 
Suppose that $f \in \Gamma$. Let optimal value $f^*$ be attained at $x^* \in C(\bx^0)$, which denotes the non-empty subset of critical points of $f$ at which Direction PPM iterates $\{\bx^k\}$ at initial $\bx^0$ can converge. Suppose that Direction PPM adopts directions of sub-gradient-relatedness and target-relatedness $\{\bar \bp_k\}$. We further suppose that $\lim\inf_k t_k = t > 0$ where $t_k$ satisfies (\ref{tlowerbound}). Denote $\bu$ as the next iterate of $\bx \in \{\bx^k\}$ with $\bu = \bx + w^* \bar \bp$.  \\
(i) $\| \bu - \bx^* \| \leq \|\bx- \bx^*\|$ holds with 
\begin{align} \label{DPPMvalue1}
2t (f(\bu) - f^* ) + (w^*)^2 \leq \|\bx - \bx^* \|^2 - \| \bu - \bx^* \|^2. 
\end{align}
(ii) The entire sequence $\{\bx^k\}$ generated by Direction PPM converges to a single critical point of $f$.  Furthermore, if $f \in \Gamma$ is differentiable, then adopting search sequence $\{\bp_k = -\nabla f(\bx^k)\}$ means that the entire sequence of Direction PPM iterates converges to a single critical point of $f$.
\end{citedthm}

\proof 
See Appendix \ref{dppmfejerproof0}.

\qed

\subsection{Rate of convergence to sub-optimal solutions} \label{convergencerate}

In the following, we examine the rate at which Direction PPM converges in terms of the number of updates to reach an $\epsilon$-suboptimal solution of $f$. 

\begin{citedlem} \label{dppmexactconvex}
Suppose that the assumption pertaining to Theorem \ref{dppmfejer} holds. 
Then, Direction PPM achieves an $\epsilon$-suboptimal solution (i.e., $f(\bx^k) -f^* \leq \epsilon$) in the order of ${\cal O}(\frac{1}{\epsilon})$. 
Furthermore, this order can be reduced to ${\cal O}({\frac{1}{\sqrt{\epsilon}}})$.
\end{citedlem}
\proof
\begin{align} 
2kt(f(\bx^k) - f^*) & \leq 
2kt(f(\bx^k) - f^*)   + \sum_{i=0}^{k-1} (w_i^*)^2  \nonumber \\ 
&\leq \sum _{i=0}^{k-1} 2t (f(\bx^{i+1}) - f(\bx^*)) + \sum_{i=0}^{k-1} (w_i^*)^2  \nonumber \\
&  \leq \sum_{i=0}^{k-1} \|\bx^{i} - \bx^{*} \| - \|\bx^{i+1} - \bx^{*} \| \nonumber \\
& \leq \|\bx^0 - \bx^{*} \| - \| \bx^k - \bx^0 \| \leq \| \bx^0 - \bx^{*}\|. \label{dppmtargetprovide3}
\end{align}
The second inequality follows from the fact that Direction PPM is a descent algorithm. The third inequality follows from (\ref{DPPMvalue1}) in Theorem \ref{dppmfejer}.
Based on (\ref{dppmtargetprovide3}), the number of iterations $k$ required to attain an $\epsilon$-suboptimal solution is in the order of $\frac{1}{\epsilon}$. This order can be accelerated to ${\cal O}(\frac{1}{\sqrt \epsilon})$ using Nesterov's acceleration method \cite{nesterov1983method}, as shown in Appendix~\ref{sec:acc}.

\qed

\section{Experiments}
\label{sec:experiments}
We present numerical results to demonstrate the efficiency of Direction PPM in comparison to existing optimization methods on convex functions. Two 2-dimensional and two higher-dimensional functions are selected as our objective functions.

\subsection{Experiment settings}
\begin{figure}[h!]
     \centering
     \begin{subfigure}[b]{0.45\textwidth}
         \centering
      \includegraphics[width=0.6\textwidth]{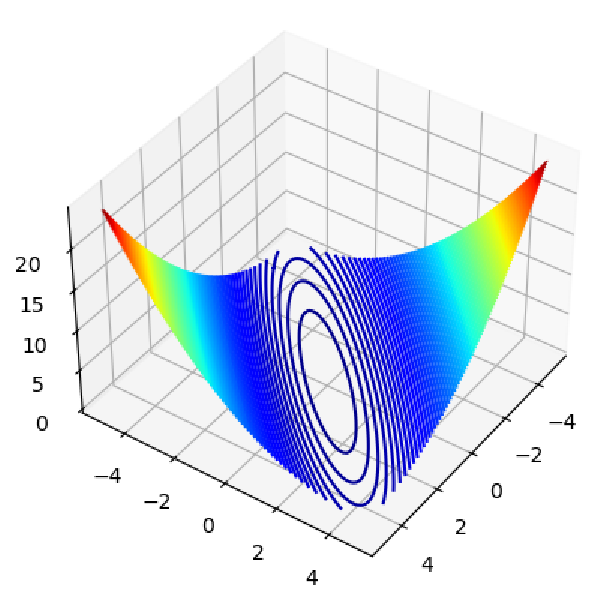}
         \label{fig:matyas}
     \end{subfigure}
         \begin{subfigure}[b]{0.45\textwidth}
         \centering
         \includegraphics[width=0.6\textwidth]{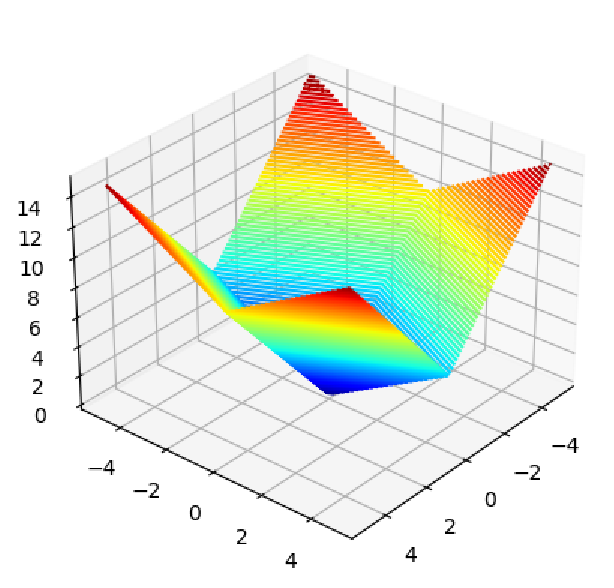}
         \label{fig:non-diff}
     \end{subfigure}
     \caption{Plots of objective functions: (left) Matyas, (right) $2|x|+|y|$}
     \label{functionplot}
\end{figure}

The first experiments involved the differentiable and non-differentiable 2-dimensional objective functions presented in Figure ~\ref{functionplot}.The differentiable objective function is a Matyas function: $f(x,y)=0.26(x^2+y^2)-0.48xy$. For Direction PPM, we set a negative gradient direction for each iteration, and use the bisection method on the directional derivative of the objective function to find the minimum for each iteration. For the baseline, we use the gradient method with backtracking line search. For the non-differentiable objective function $f(x,y)=2|x|+|y|$, the baseline algorithm is based on the sub-gradient method, in which the sub-gradient at the current iterate  is derived from the average of the gradients at points sampled in a neighborhood of the iterate.
For Direction PPM, the search direction is set the same as the above sub-gradient at each iterate, and the optimal point of scalar optimization at each iteration is derived using the Golden Search method. 


The next two examples are drawn from real-world applications. The first one is related to compressed sensing. We consider $f(\bx)=\|A \bx-\bb \|_1 + 10 \| \bx \|_1$, where $A\in\mathbb{R}^{10\times50}$ and $\bb\in\mathbb{R}^{10}$ are fixed, and $\bx\in\mathbb{R}^{50}$ is the variable that is optimized.  Each entry of $A$ is generated from independently and identically distributed (i.i.d.) Gaussian  distribution $\mathcal N(0, \frac{1}{50})$, and $\bb$ is the product of $A$ and a 50-dimensional vector with only 5 nonzero terms, where the coefficients are generated uniformly from $[-0.5, 0.5]$. In this example, we simulate the procedure of recovering the original sparse vector. The second example is a classification problem with logistic loss and a $L_1$ regulation term. 
Given input and label pairs $\{(\bx,y)^{(i)}, \bx^{(i)} \in\mathbb{R}^{10},y^{(i)} \in\{0,1\}\}$ where $\bx^{(i)}$ is uniformly generated from $[-5, 5]^{10}$, and $y^{(i)}$ is generated from a Bernoulli distribution with equal probability, we aim to find the optimal $\bw\in\mathbb{R}^{10}$ that minimizes $f(\bw)=\frac{1}{100}\sum_{i=1}^{100}\log(1 + e^{-(\bw^{\top}\bx)y})+\lambda\|\bw\|_1$. Here, $\lambda$ is the coefficient of regulation term, and we set $\lambda=50$ for the following experiments. For both problems, we introduce the momentum method \cite{QIAN1999145} to select the search direction for Direction PPM. This means that the direction is chosen to be a linear combination of the direction vector in the previous step and the sub-gradient at the current iterate, in which the weights for the combination are determined experimentally.

\subsection{Results}

\begin{figure}[h!]
     \centering
     \begin{subfigure}[b]{0.45\textwidth}
         \centering
         \includegraphics[width=0.8\textwidth]{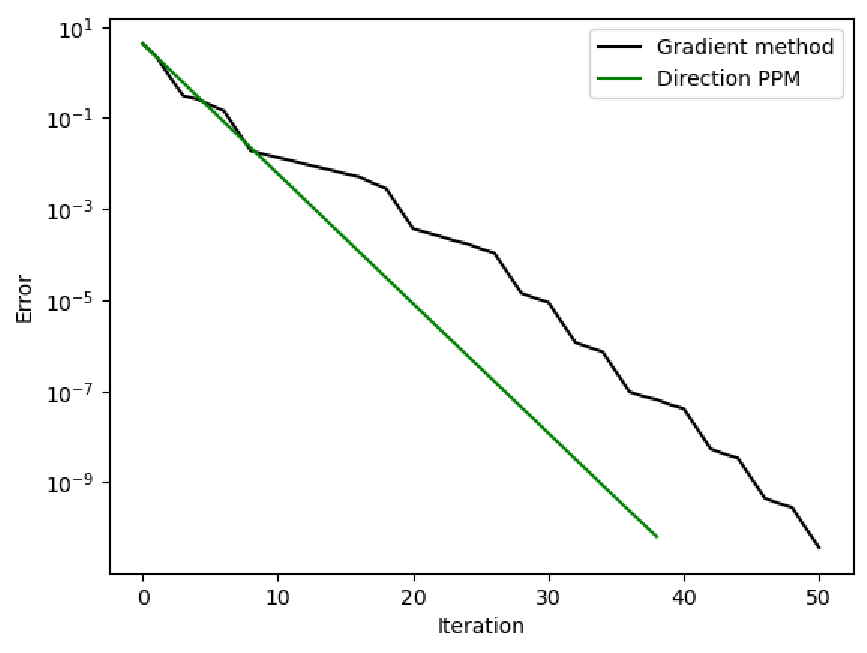}
         \label{fig:matyas_result}
     \end{subfigure}
     \begin{subfigure}[b]{0.45\textwidth}
         \centering
         \includegraphics[width=0.8\textwidth]{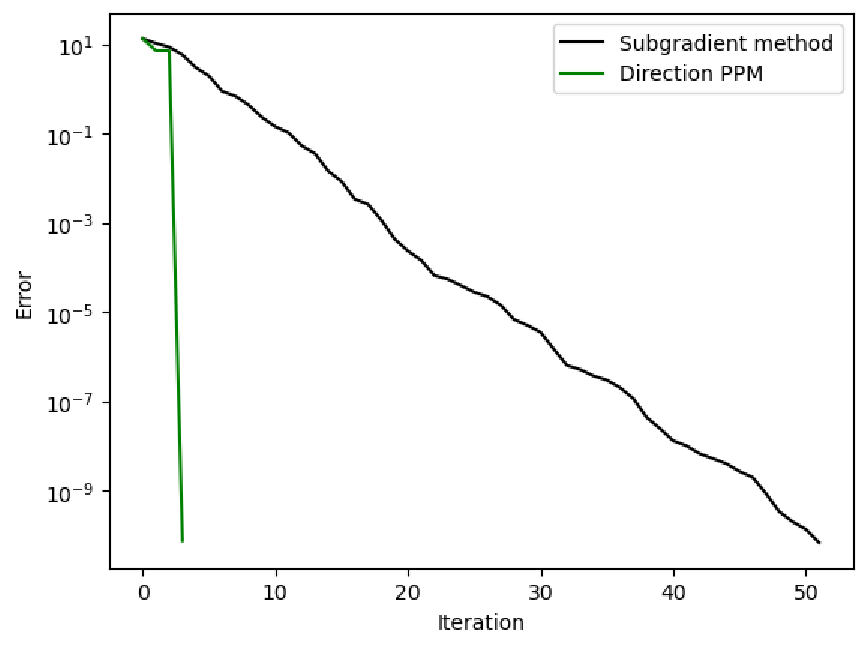}
         \label{fig:non-diff_result}
     \end{subfigure}
     \caption{Convergence of algorithms for (left) Matyas and (right) $2|x|+|y|$ where vertical axis is true error ($f(x,y)-f^*$ where $f^*$ is the optimal value) and horizontal axis is iteration number with parameter $t$ in Direction PPM set at $1000$.}
\end{figure}
\begin{table}[h!]
\centering
 \begin{tabular}{||c c c||}
 \hline
  & Sub-gradient  & Direction PPM  \\ 
 \hline\hline
 Matyas & $1.875$   & $0.781$ \\ 
 $2|x|+|y|$ &  $2.656$ & $0.469$   \\
 \hline
 \end{tabular}
 \captionof{table}{Comparison of run time (ms) for different functions and algorithms}
  \label{table:time}
\end{table}

\begin{figure}[h!]
     \centering
     \begin{subfigure}[b]{0.45\textwidth}
         \centering
         \includegraphics[width=0.8\textwidth]{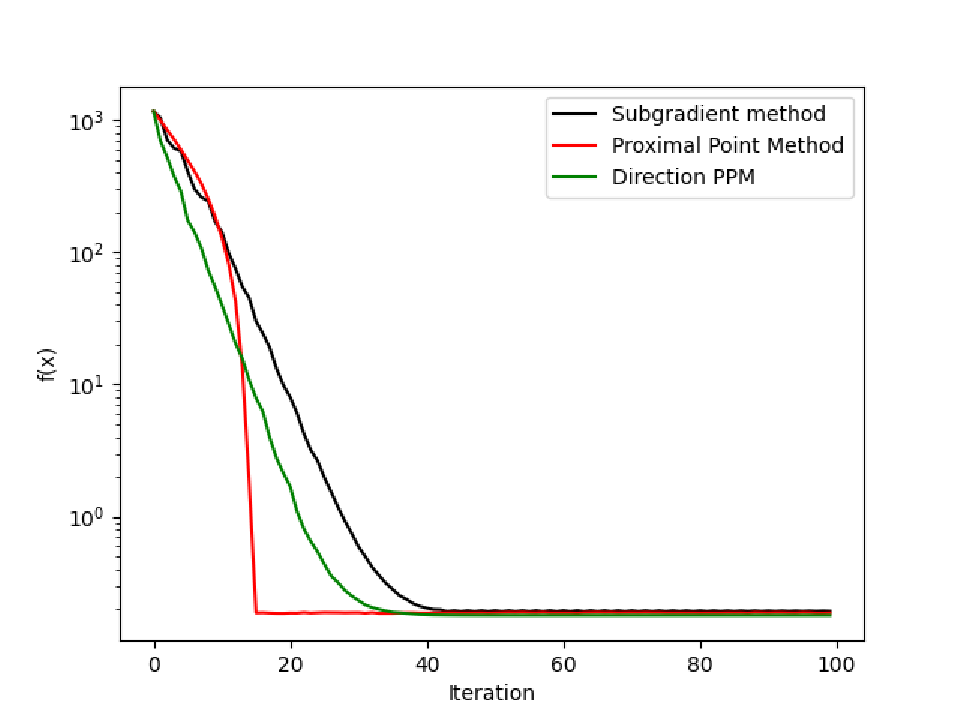}
         \label{fig:11norm_result}
     \end{subfigure}
     \begin{subfigure}[b]{0.45\textwidth}
         \centering
         \includegraphics[width=0.8\textwidth]{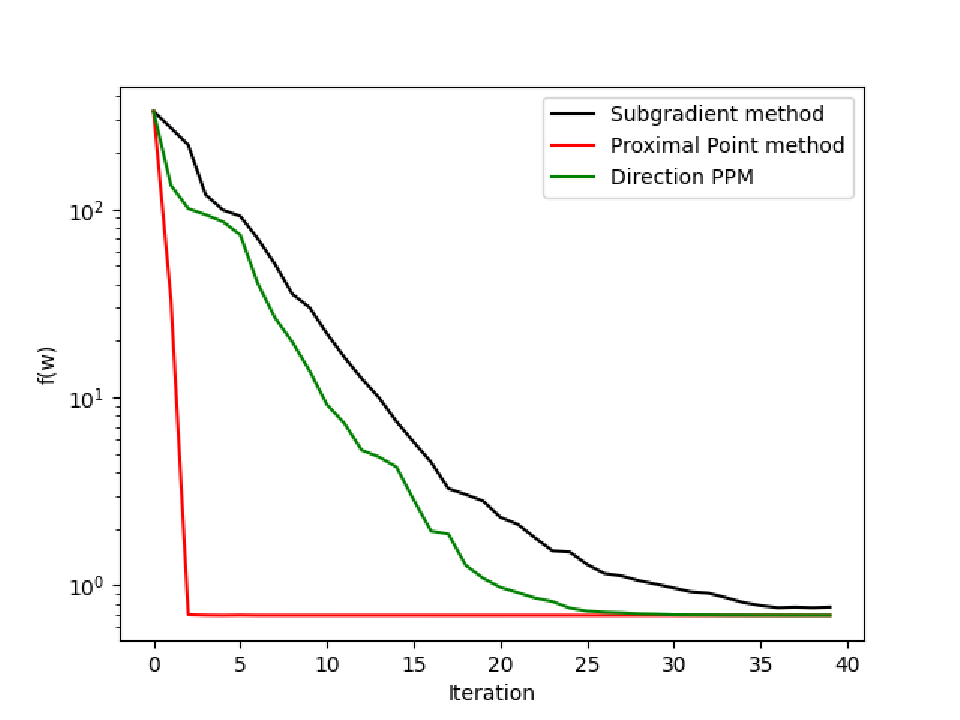}
         \label{fig:classification_result}
     \end{subfigure}
     \caption{Convergence of algorithms with artificial datasets where vertical axis is $f(\bx)$ and horizontal axis is iteration number with parameter $t$ in Direction PPM set at $1000$: (left) compressed sensing, (right) logistic regression.  }
     \label{fig:error}
\end{figure}

\begin{table}[h!]
\centering
 \begin{tabular}{||c c c c||}
 \hline
  & Sub-gradient & PPM & Direction PPM \\
 \hline\hline
 Compressed sensing & $0.194$ &  $0.186$ & $0.184$ \\
 Logistic regression & $0.767$ & $0.693$ & $0.695$\\ 
 \hline
 \end{tabular}
 \captionof{table}{Terminated value for different functions and algorithms}
 \label{table:error}
\end{table}

\begin{figure}[h!]
     \centering
     \begin{subfigure}[b]{0.45\textwidth}
         \centering
         \includegraphics[width=0.8\textwidth]{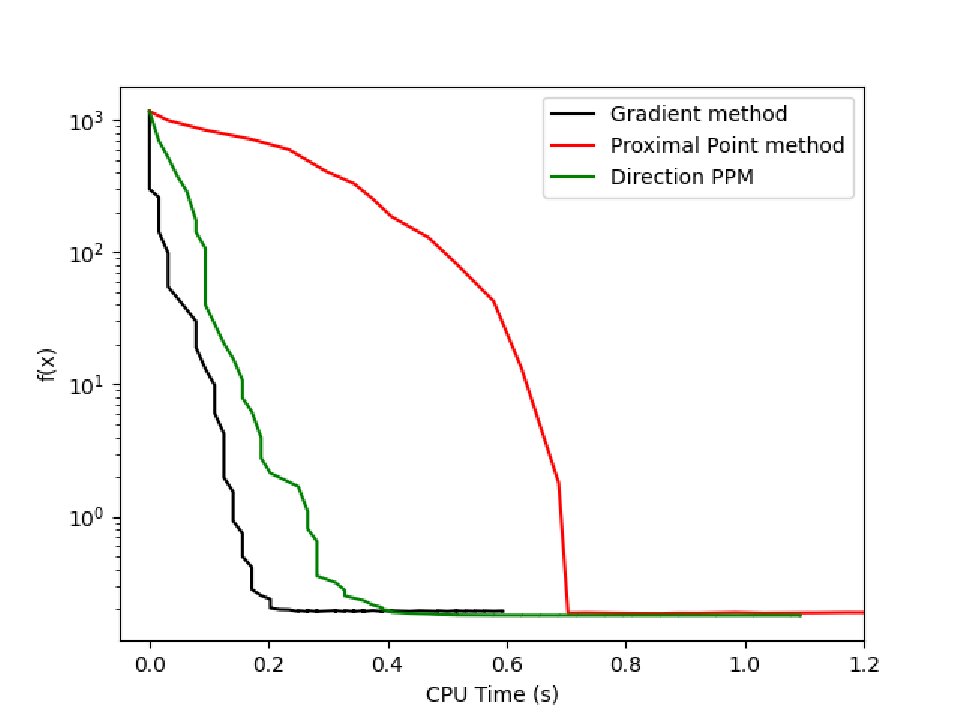}
         \label{fig:11norm_result}
     \end{subfigure}
     \begin{subfigure}[b]{0.45\textwidth}
         \centering
         \includegraphics[width=0.8\textwidth]{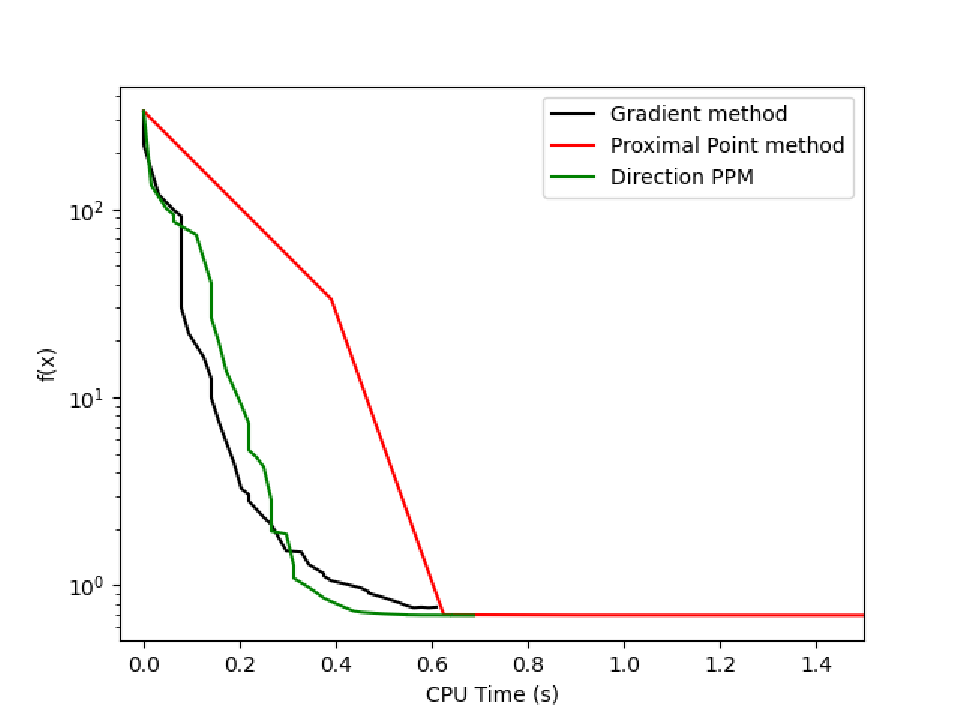}
         \label{fig:classification_result}
     \end{subfigure}
     \caption{Run time of experiments with artificial datasets: (left) compressed sensing, (right) logistic regression}
     \label{time}
\end{figure}

We present our experimental results in the following graphs. All the experiments were executed on an Intel(R) Core(TM) i5-8265U CPU. Note that for the 2-dimensional functions, because the optimal values are zeros, we could calculate the true errors for those two algorithms. For the compressed sensing and logistic regression objective functions, the optimal values are unknown. Therefore we can only compare the relative performance of different algorithms, i.e., the terminated value of functions on the corresponding sub-optimal points derived by algorithms. Furthermore, we present the run time for each algorithm to illustrate that the proposed algorithm is not computationally expensive. For the 2-dimensional cases, we measure how long the function takes to converge to $10^{-10}$. As it is difficult to determine the convergence time for the high-dimensional cases, we provide function value versus CPU time plots to visualize the speed at the algorithms converge.

The 2-dimensional examples show the potential of Direction PPM when applied to difficult optimization problems. Because two eigenvalues of the Heissian matrix in the Matyas function differ significantly from one another, the gradient-based algorithm usually suffers from choosing step sizes to avoid zig-zag paths. The non-smoothness of function $2|x|+|y|$ often confuses solvers and slows convergence in a theoretical aspect. The results presented in Table ~\ref{table:time} show that Direction PPM takes shorter time to converge in the same level as sub-gradient method does, which means that Direction PPM overcomes the common obstacles of optimization problems.

The higher-dimensional optimization problems are more challenging to solve. PPM is used as baseline algorithm. For each iteration, we use gradient descent to solve the subproblem with respect to the proximal operator (due to the fact that Moreau envelope is differentiable and its gradient is $\frac{1}{t}$ Lipschitz continuous \cite{BC2011}). For the gradient descent procedure, we set the step size as $1/i^{1.5}$ for the $i$-th descent, and we run $150$ descents for compressed sensing and $400$ for logistic regression. For the proximal operator, we set parameter $t=1000$. As Figure ~\ref{fig:error} and Table ~\ref{table:error} show, this algorithm solves these two problems well. However, its computational cost is high, as every iteration contains an optimization problem of almost the same complexity as the original problem. While the sub-gradient method is simpler and more efficient, the final convergence value is significantly worse than that of PPM. The proposed Direction PPM preserves the advantages of both methods. As shown in Figure ~\ref{time}, for the compressed sensing example, PPM takes around 0.7 seconds to converge, while Direction PPM only takes 0.4 seconds. Similarly, in the logistic regression example, PPM method take around 0.6 seconds to converge, while Direction PPM only takes 0.45 seconds. These two examples indicate that Direction PPM is around 30\% to 40\% faster than PPM. 

\section{Conclusions} \label{cons}

The use of proximal point operators for optimization can be computationally expensive when the dimensionality of a function (i.e., the number of variables) is high. We reduce the cost of calculating the proximal point operators by splitting the search into two separate tasks: determining the direction from one iterate to the next, and searching for a suitable step-size as a scalar optimization problem along the chosen direction.
We show that under certain conditions, the proposed Direction PPM achieves the optimal points at a guaranteed rate. 
We also demonstrate numerically that the computational cost of  Direction PPM is significantly less than that of  PPM, with a similar level of performance. The key to this improvement lies in reducing the dimension of each subproblem. Thus, selection of the descent direction is crucial. In future research, we plan to focus on direction selection that is flexible enough to cope with large families of functions.

%

\bibliographystyle{ieeetr}
\bibliography{DeepRef}

\begin{thebibliography}{10}

\bibitem{bolte2014proximal}
J.~Bolte, S.~Sabach, and M.~Teboulle, ``Proximal alternating linearized
  minimization for nonconvex and nonsmooth problems,'' {\em Mathematical
  Programming}, vol.~146, no.~1, pp.~459--494, 2014.

\bibitem{becker2011nesta}
S.~Becker, J.~Bobin, and E.~J. Cand{\`e}s, ``Nesta: A fast and accurate
  first-order method for sparse recovery,'' {\em SIAM Journal on Imaging
  Sciences}, vol.~4, no.~1, pp.~1--39, 2011.

\bibitem{BoydADMMsurvey}
S.~Boyd, N.~Parikh, E.~Chu, B.~Peleato, and J.~Eckstein, ``Distributed
  optimization and statistical learning via the alternating direction method of
  multipliers,'' {\em Found. Trends Mach. Learn.}, vol.~3, no.~1, pp.~1--122,
  2011.

\bibitem{zhang2018ista}
J.~Zhang and B.~Ghanem, ``Ista-net: Interpretable optimization-inspired deep
  network for image compressive sensing,'' in {\em Proceedings of the IEEE
  conference on computer vision and pattern recognition}, pp.~1828--1837, 2018.

\bibitem{parikh2014proximal}
N.~Parikh, S.~Boyd, {\em et~al.}, ``Proximal algorithms,'' {\em Foundations and
  trends{\textregistered} in Optimization}, vol.~1, no.~3, pp.~127--239, 2014.

\bibitem{parikh2014block}
N.~Parikh and S.~Boyd, ``Block splitting for distributed optimization,'' {\em
  Mathematical Programming Computation}, vol.~6, no.~1, pp.~77--102, 2014.

\bibitem{asl2020analysis}
A.~Asl and M.~L. Overton, ``Analysis of the gradient method with an
  armijo--wolfe line search on a class of non-smooth convex functions,'' {\em
  Optimization methods and software}, vol.~35, no.~2, pp.~223--242, 2020.

\bibitem{fletcher2005barzilai}
R.~Fletcher, ``On the barzilai-borwein method,'' in {\em Optimization and
  control with applications}, pp.~235--256, Springer, 2005.

\bibitem{nocedal2006numerical}
J.~Nocedal and S.~Wright, {\em Numerical optimization}.
\newblock Springer Science \& Business Media, 2006.

\bibitem{Ber08}
D.~P. Bertsekas, ``Nonlinear programming: 3rd,'' {\em Athena Scientific
  Optimization and Computations Series 4}, vol.~4, 2008.

\bibitem{barzilai1988two}
J.~Barzilai and J.~M. Borwein, ``Two-point step size gradient methods,'' {\em
  IMA journal of numerical analysis}, vol.~8, no.~1, pp.~141--148, 1988.

\bibitem{burdakov2019stabilized}
O.~Burdakov, Y.-H. Dai, and N.~Huang, ``Stabilized barzilai-borwein method,''
  {\em Journal of Computational Mathematics}, pp.~916--936, 2019.

\bibitem{bottou2018optimization}
L.~Bottou, F.~E. Curtis, and J.~Nocedal, ``Optimization methods for large-scale
  machine learning,'' {\em Siam Review}, vol.~60, no.~2, pp.~223--311, 2018.

\bibitem{kinderlehrer2000introduction}
D.~Kinderlehrer and G.~Stampacchia, {\em An introduction to variational
  inequalities and their applications}.
\newblock SIAM, 2000.

\bibitem{press1988numerical}
W.~H. Press, S.~A. Teukolsky, W.~T. Vetterling, and B.~P. Flannery, ``Numerical
  recipes in c,'' 1988.

\bibitem{avriel1968golden}
M.~Avriel and D.~J. Wilde, ``Golden block search for the maximum of unimodal
  functions,'' {\em Management Science}, vol.~14, no.~5, pp.~307--319, 1968.

\bibitem{xia1992finding}
Z.-q. Xia, ``Finding subgradients or descent directions of convex functions by
  external polyhedral approximation of subdifferentials,'' {\em Optimization
  Methods and Software}, vol.~1, no.~3, pp.~253--264, 1992.

\bibitem{bauschke2011convex}
H.~H. Bauschke, P.~L. Combettes, {\em et~al.}, {\em Convex analysis and
  monotone operator theory in Hilbert spaces}, vol.~408.
\newblock Springer, 2011.

\bibitem{nesterov1983method}
Y.~Nesterov, ``A method for unconstrained convex minimization problem with the
  rate of convergence o (1/k\^{} 2),'' in {\em Doklady an ussr}, vol.~269,
  pp.~543--547, 1983.

\bibitem{QIAN1999145}
N.~Qian, ``On the momentum term in gradient descent learning algorithms,'' {\em
  Neural Networks}, vol.~12, no.~1, pp.~145--151, 1999.

\bibitem{BC2011}
H.~H. Bauschke and P.~L. Combettes, {\em Convex Analysis and Monotone Operator
  Theory in Hilbert Spaces}.
\newblock Springer Verlag, 2011.

\end{thebibliography}

\appendix

\section{Proof of Theorem \ref{dppmfejer}} \label{dppmfejerproof0}

(i) Proposition \ref{dppmlimitpoints} allows us to assign $f(\bx^*) = f^*$ for any $x^* \in C(\bx^0)$. 
Let $\bx^* \in C(\bx^0)$ and let $\bu = \bx- t (\bs(\bu)^\top \bar \bp) \bar \bp$ where $\bs(\bu) \in \partial f(\bu)$.
Then, we have
\begin{align} \label{DPPMdistance}
\| \bu - \bx^* \|^2 & = \| \bx - t (\bs(\bu)^\top \bar \bp) \bar \bp- \bx^* \|^2  \nonumber \\
& = \|\bx - \bx^* \|^2 - 2t (\bs(\bu)^\top \bar \bp) \bar \bp^\top (\bx - \bx^*) + t^2 |\bs(\bu)^\top \bar \bp |^2.
\end{align}
Using $f$ is convex, 
we obtain 
\begin{align}
f(\bz) \geq f(\bu) + \bs(\bu)^\top (\bz - \bu) 
 = f(\bu)  + \bs(\bu)^\top  (\bz - \bx + t(\bs(\bu)^\top \bar \bp) \bar \bp).
\end{align}
This equation can be re-written as  
\begin{align} \label{DPPMmain0}
f(\bu)  \leq  f(\bz)  + \bs(\bu)^\top (\bx  - \bz) - t | \bs(\bu)^\top \bar \bp|^2.
\end{align}
Applying (\ref{gradientdecomp}) to $\bs(\bu)^\top (\bx  - \bz)$, we obtain  
$\bs(\bu)^\top (\bx - \bz) = (\alpha \bar \bp + \beta \bv)^\top (\bx - \bz)$.
This result is then substituted into (\ref{DPPMmain0}) to yield
\begin{align} \label{fejer1}
f(\bu) \leq f(\bz) + (\alpha \bar \bp + \beta \bv)^\top (\bx - \bz)   - t | \bs(\bu)^\top \bar \bp|^2.
\end{align}
From (\ref{gradientdecomp}), we obtain
$\alpha = \bar \bp^\top \bs(\bu) - \beta \bar \bp^\top \bv$.
Substituting this $\alpha$ into (\ref{fejer1}) yields
\begin{align} \label{fejer11}
f(\bu)  \leq  f(\bz) + (\bar \bp^\top \bs(\bu) - \beta \bar \bp^\top \bv) \bar \bp + \beta \bv)^\top (\bx - \bz)   - t | \bs(\bu)^\top \bar \bp|^2.
\end{align}
We can now substitute $\bx^*$ for $\bz$ to obtain  
\begin{align} \label{fejer2}
f(\bu) & \leq  f^* +  (\bs(\bu)^\top \bar \bp) \bar \bp^\top + \beta \bv^\top - (\beta  \bar \bp^\top \bar \bv) \bar \bp^\top)(\bx - \bx^*) - t | \bs(\bu)^\top \bar \bp|^2.
\end{align}
Multiplying both sides of (\ref{fejer2}) by $2t$, and replacing $ 2t (\bs(\bu)^\top \bar \bp) \bar \bp^\top (\bx - \bx^*)$ in accordance with (\ref{DPPMdistance}), we obtain 
\begin{align} \label{DPPMvalue0}
2t (f(\bu) - f^*) & \leq  \|\bx - \bx^* \|^2 - \| \bu - \bx^* \|^2 - t^2 |\bs(\bu)^\top \bar \bp |^2 \nonumber \\
& + 2 t\beta \bv^\top (\bx - \bx^*) - 2t \beta (\bar \bp^\top \bv) \bar \bp^\top (\bx - \bx^*)
\end{align}
Since $\bx - \bx^* = -\bar \bp^* \|\bx - \bx^*\|$ and (\ref{equalproj}) (under the assumption of target-relatedness of $\bar \bp$), we can obtain
\begin{align} \label{DPPMvalue4}
2t \beta  (\bv^\top  - (\bar \bp^\top \bv) \bar \bp^\top) (\bx - \bx^*) &  = - 2 \|\bx - \bx^*\|t \beta  (\bv^\top  - (\bar \bp^\top \bv) \bar \bp^\top)\bar \bp^* \nonumber \\
 & = - 2 \|\bx - \bx^*\|t \beta   (\bv^\top \bar \bp^* - (\bar \bp^\top \bv)^2 ) \nonumber \\
 & =  2 \|\bx - \bx^*\|t \beta   (-\| \bv \|^2 + (\bar \bp^\top \bv)^2 ) 
\end{align}
Because $-\| \bv \|^2 + (\bar \bp^\top \bv)^2  \leq 0$, if $\beta \geq 0$, then (\ref{DPPMvalue4}) is smaller than zero; hence, (\ref{DPPMvalue0}) becomes
\begin{align} \label{DPPMvalue2}
2t (f(\bu) - f^*) & \leq  \|\bx - \bx^* \|^2 - \| \bu - \bx^* \|^2 - t^2 |\bs(\bu)^\top \bar \bp |^2.
\end{align}
Note that under the assumption of target-relatedness, we obtain $\bar \bp^\top \bar \bv > 0$ from which to obtain $\beta \geq 0$ in accordance with Eq. (\ref{cosine}).
Because $w^* = t |\bs(\bu)^\top \bar \bp|$ (Lemma \ref{dppmunique}), (\ref{DPPMvalue2}) can be expressed as 
\begin{align*} 
0 \leq 2t (f(\bu) - f^* ) + (w^*)^2 \leq \|\bx - \bx^* \|^2 - \| \bu - \bx^* \|^2. 
\end{align*}

(ii) Since the DPPM-directions satisfy the assumption of sub-gradient-relatedness, we can follow Proposition \ref{dppmlimitpoints} and suppose that sub-sequences $\{\bx^{y(k)}\}_k$ and $\{\bx^{z(k)}\}_k$ respectively converge to critical points $\by$ and $\bz$. In accordance with the assumption of target-relatedness and (i), sequences $\{\|\bx^k - \by\| \}$ and  $\{\|\bx^k - \bz\| \}$ converge. The following derivations can be found in \cite[Lemma 2.39]{bauschke2011convex}. From
\begin{align*} 
2 (\bx^k)^\top (\bz-\by) = \|\bx^k - \by \|^2 - \| \bx^k - \bz\|^2  + \|\bz\|^2 - \|\by\|^2 \text{ for all $k$,}
\end{align*}
we can deduce  that $\{(\bx^k)^\top (\bz-\by)\}$ also converges, say $\{(\bx^k)^\top (\bz-\by)\} \rightarrow l$. Proceeding to the limits along $\{\bx^{y(n)}\}$ and $\{\bx^{z(n)}\}$ respectively yields
\begin{align} \label{entire1}
l = \by^\top (\bz - \by)
\end{align}
and 
\begin{align} \label{entire2}
l = \bz^\top (\bz - \by).
\end{align}
From (\ref{entire1}) and (\ref{entire2}), $\| \by - \bz \|^2 = 0$; therefore, $\by = \bz$.

%
%
%
%
%
%
%

\section{Direction PPM Acceleration} \label{sec:acc}

For the sake of convenience, Eqs. (\ref{dppmdescent})  and (\ref{DPPMvalue1}) are re-stated as follows:
\begin{align}  \label{dppmaccer1}
f(\bu) \leq f(\bx) - t | \nabla f(\bu)^\top \bar \bp|^2, 
\end{align}
\begin{align}  \label{dppmaccer0}
t( f(\bu) -f^*)   \leq \|\bx - \bx^* \|^2 - \| \bu - \bx^* \|^2.
\end{align} 

Let $z =  | \nabla f(\bu)^\top \bar \bp|$. Multiply (\ref{dppmaccer1}) by $t(1-\theta)$ where $\theta \in [0, 1]$, and multiply (\ref{dppmaccer0}) by $\theta$. Add the resultants together to yield the following:
\begin{align*}
t f(\bu) & \leq (1-\theta) t f(\bx) + t \theta f^* - t^2 (1-\theta) z^2 + \theta (\|\bx - \bx^* \|^2 - \| \bu - \bx^* \|^2) \\
& \leq (1-\theta) t f(\bx) + t \theta f^* - t^2 (1-\theta) z^2 - 2 tz \theta \bar \bp^\top (\bx - \bx^*) - \theta t^2 z^2.
\end{align*}
The last equality is derived with $\bu= \bx - t (\nabla f(\bu)^\top \bar \bp) \bar \bp = \bx + t z \bar \bp$ in accordance with $\nabla f(\bu)^\top \bar \bp \leq 0$.
Following the inequality, we obtain
\begin{align*}
f(\bu) & \leq (1- \theta) f(\bx) + \theta f^* - 2\theta z \bar \bp^\top (\bx - \bx^*) - t z^2 \\
& = (1- \theta) f(\bx) + \theta f^* + \frac{\theta^2}{t} (\|\bx - \bx^* \|^2 - \| \bx + \frac{tz \bar \bp}{\theta} - \bx^* \|^2).
\end{align*}
Hence, 
\begin{align} \label{dppmaccer2}
f(\bu) - f^* \leq (1-\theta) f(\bx) - (1-\theta) f^* + \frac{\theta^2}{t}  (\|\bx - \bx^* \|^2 - \| \bx + \frac{tz \bar \bp}{\theta} - \bx^* \|^2).
\end{align}
Letting $\bv = \bx + \frac{tz \bar \bp}{\theta}$ and recall that $
\bu = \bx + tz \bar \bp$,  we obtain  
\begin{align} \label{vextra}
\bv = \bx + \frac{1}{\theta} (\bu - \bx) = (1 - \frac{1}{\theta}) \bx + \frac{1}{\theta} \bu.
\end{align}

Dividing both sides of (\ref{dppmaccer2}) by $\theta^2$ and re-arranging the terms yields the following:
\begin{align} \label{dppmaccer2}
\frac{1}{\theta^2} (f(\bu) - f^*) + \frac{1}{t} \|\bv - \bx^*\|^2 \leq \frac{1-\theta}{\theta^2}[f(\bx) - f^*] + \frac{1}{t} \|\bx - \bx^*\|^2.
\end{align}
If we retrieve the iteration index of the above and replace $\bu, \bv, \bx, \theta$ respectively with $\bx^k$, $\bv^k$, $\bv^{k-1}$, and $\theta_k$, Eq. (\ref{dppmaccer2}) can be expressed as follows:
\begin{align} \label{dppmaccer3}
\frac{1}{\theta_k^2} (f(\bx^k) - f^*) + \frac{1}{t} \|\bv^k - \bx^*\|^2 \leq \frac{1-\theta_k}{\theta_k^2}[f(\bv^{k-1}) - f^*] + \frac{1}{t} \|\bv^{k-1} - \bx^*\|^2.
\end{align}
Letting 
\begin{align}
\theta_i = \frac{2}{i+1},
\end{align}
we have the desired properties: $\theta_1 = 1$, $\theta_i \in (0, 1)$, and 
\begin{align} \label{dppmthetacon}
\frac{1-\theta_i}{\theta_i^2} \leq \frac{1}{\theta_{i-1}^2}.
\end{align}
Repeatedly using (\ref{dppmaccer3}) and (\ref{dppmthetacon}), we obtain
\begin{align} \label{dppmaccer4}
\frac{1}{\theta_k^2} (f(\bx^k) - f^*) + \frac{1}{t} \|\bv^k - \bx^*\|^2 & \leq \frac{1}{\theta_{k-1}^2}[f(\bv^{k-1}) - f^*] + \frac{1}{t} \|\bv^{k-1} - \bx^*\|^2 \nonumber \\
& \leq \frac{1-\theta_{k-1}}{\theta_{k-1}^2}[f(\bv^{k-2}) - f^*] + \frac{1}{t} \|\bv^{k-2} - \bx^*\|^2 \nonumber\\
& \leq \frac{1}{\theta_{k-2}^2}[f(\bv^{k-2}) - f^*] + \frac{1}{t} \|\bv^{k-2} - \bx^*\|^2 \nonumber\\
& \vdots \nonumber\\
& \leq \frac{1}{\theta_{1}^2}[f(\bv^{1}) - f^*] + \frac{1}{t} \|\bv^{1} - \bx^*\|^2 \nonumber\\
& \leq \frac{1- \theta_1}{\theta_{1}^2}[f(\bv^{0}) - f^*] + \frac{1}{t} \|\bv^{0} - \bx^*\|^2 \nonumber \\
& \leq \frac{1}{t} \|\bv^{0} - \bx^*\|^2.
\end{align}
Deduced from (\ref{dppmaccer4}) and $\bv^0 = \bx^0$,  
\begin{align} \label{dppmaccer5}
f(\bx^k) - f^* \leq \frac{\theta_k^2}{t} \|\bx^{0} - \bx^*\|^2.
\end{align}
The order of $k$ to achieve $\epsilon$-suboptimal solution of $f^*$ is $\frac{1}{\sqrt{\epsilon}}$.

\qed

Eq. (\ref{vextra}) indicates that $\bv^k$ is an extrapolation of $\bx^k$ and $\bu^k$, where parameter $\theta_k = \frac{2}{k+1}$ where the initial setting at $\bx^0 = \bv^0$.

\end{document}